\documentclass{amsart}
\usepackage{amsmath,amsthm,amsfonts,amssymb}
\date{\today}

\newcommand{\C}{\mathcal C}
\newcommand{\IN}{\mathbb N}

\newcommand{\w}{\omega}
\newcommand{\cc}{\mathfrak c}
\newcommand{\Ra}{\Rightarrow}

\newtheorem{theorem}{Theorem}

\newtheorem{proposition}{Proposition}
\newtheorem{problem}{Problem}
\theoremstyle{definition}
\newtheorem{remark}{Remark}

\title[The Rees-Suschkewitsch Theorem]
{The Rees-Suschkewitsch Theorem for simple topological semigroups}

\author{Taras Banakh}
\address{Department of
Mathematics, Lviv National University, Universytetska 1, Lviv, 79000,
Ukraine} \email{tbanakh@yahoo.com}

\author{Svetlana~Dimitrova}
\address{National Technical University "Kharkov Polytechnical Institute", Frunze 21, Kharkiv, 61002, Ukraine}

\author{Oleg~Gutik}
\address{Department of Mechanics and Mathematics, Ivan Franko Lviv
National University, Universytetska 1, Lviv, 79000, Ukraine}
\email{o\_\,gutik@franko.lviv.ua}

\begin{document}



\begin{abstract} 
%
%
%
%
We detect topological semigroups that are topological paragroups, i.e., are isomorphic to a Rees product $[X\times H\times Y]_\sigma$ of a topological group $H$ over topological spaces $X,Y$ with a continuous sandwich function $\sigma:Y\times X\to H$. We prove that a simple topological semigroup $S$ is a topological paragroup if one of the following conditions is satisfied: (1) $S$ is completely simple and the maximal subgroups of $S$ are topological groups, (2) $S$ contains an idempotent and the square $S\times S$ is countably compact or pseudocompact, (3) $S$ is sequentially compact or the power  $S^{2^{\mathfrak c}}$ is  countably compact.
The last item generalizes an old Wallace's result saying that each simple compact topological semigroup is a topological paragroup.
\end{abstract}

\maketitle

This paper was motivated by the classical Rees-Suschkewitsch
Theorem that describes the algebraic structure of completely
simple semigroups and the topological versions   of this theorem
proved for compact topological semigroups by
Wallace~\cite{Wallace}, for compact semitopological semigroups by
Ruppert~\cite{Rup}, and for sequential countably compact
topological semigroups by Gutik, Pagon and Repov\v{s}~\cite{GPR}.
All topological semigroups considered in this paper are Hausdorff.

We recall that a semigroup $S$ is {\em simple} if $S$ contains no
proper two-sided ideal. A simple semigroup $S$ is called {\em
completely simple} if the set  $E=\{e\in S : ee=e\}$ of
idempotents of $S$ contains a {\em primitive idempotent}, that is,
a  minimal idempotent with respect to the partial order $e\le f$
on $E$ defined by $ef=fe=e$. In this case all the idempotents are
primitive and $H_e=eSe$ is a group for every $e\in E$, see
\cite[Section~2.7, Ex.~6(b)]{CP} or \cite{Wallace1957}.

Let us observe that each group is a completely simple semigroup.
Less trivial examples of such semigroups appear as minimal ideals
in compact right-topological semigroups, see
\cite[Theorem~I.3.13]{Rup}, \cite[Theorem~2.9]{HS} and
\cite{Ruppert}. A generic example of a completely simple semigroup
can be constructed as follows. Take any group $H$ and a function
$\sigma:Y\times X\to H$ defined on the product of two sets. This
function $\sigma$ induces the semigroup operation
$$(x,h,y)\cdot(x',h',y')=(x,h\sigma(y,x')h',y')$$ on the product
$X\times H\times Y$ turning it into a completely simple semigroup,
called the \emph{Rees product of $H$ over $X$ and $Y$ relative to
the sandwich map} $\sigma$~\cite{HofmannMostert} or a
\emph{paragroup}~\cite{Rup} and denoted by $[X,H,Y]_\sigma$.

The Rees-Suschkewitsch Structure Theorem \cite{Rees1940} says that
the converse is also true: each completely simple semigroup $S$ is
isomorphic to the paragroup $[X_e,H_e,Y_e]_\sigma$ where $e$ is
any idempotent of $S$, $H_e=eSe$ is the maximal subgroup of $S$
containing $e$, $X_e=Se\cap E$, $Y_e=eS\cap E$, and  the sandwich
function $\sigma:Y_e\times X_e\to H_e$ is defined by
$\sigma(y,x)=yx$. In fact, the map $$R:[X_e,H_e,Y_e]_\sigma\to
S,\;\; R:(x,h,y)\mapsto xhy,$$ is an isomorphism called the {\em
Rees isomorphism}. Its inverse $R^{-1}:S\to [X_e,H_e,Y_e]_\sigma$
is defined by the formula
$$R^{-1}(s)= (s(ese)^{-1},ese,(ese)^{-1}s).$$

Now assume that $S$ is a topological semigroup (i.e., a
topological space endowed with a continuous semigroup operation).
In this case the spaces $X_e=Se\cap E$, $H_e=eSe$, and $Y_e=eS\cap
E$ carry the induced topologies while the paragroup
$[X_e,H_e,Y_e]_\sigma$ carries the product topology making the
semigroup operation continuous, i.e., $[X_e,H_e,Y_e]_\sigma$ is a
topological semigroup. Let us observe that the maximal subgroup
$H_e$ of $S$ is a paratopological group, which means that the
group multiplication on $H_e$ is jointly continuous. If, in
addition, the inversion map $x\mapsto x^{-1}$ is continuous on
$H_e$, then $H_e$ is a {\em topological group}.

Looking at the Rees isomorphism
$$R:[X_e,H_e,Y_e]_\sigma\to S,\;\;R:(x,h,y)\mapsto xhy,$$ we see
that it is continuous while its inverse
$$R^{-1}:s\mapsto (s(ese)^{-1},ese,(ese)^{-1}s)$$ is continuous if
the paratopological group $H_e$ is a topological group. In this
case the topological semigroup $[X_e,H_e,Y_e]_\sigma$ is called a
\emph{topological paragroup}~\cite{HofmannMostert}.

More precisely, by {\em a topological paragroup} we understand a
topological semigroup that is topologically isomorphic to the Rees
product $[X,H,Y]_\sigma$ where $H$ is a topological group and
$\sigma:X\times Y\to H$ is a continuous function defined on the
product of two topological spaces.

In such a way we have obtained the following topological
Rees-Suschkewitsch structure theorem.

\begin{theorem}\label{t1}
A topological semigroup $S$ is a topological paragroup if and only
if $S$ is completely simple and each maximal subgroup $H_e$ of $S$
is a topological group.
\end{theorem}

There is a simple algebraic characterization of completely simple
semigroups, see  \cite[Theorem 2.54]{CP} and \cite{Andersen}.

\begin{theorem}[Andersen]\label{an}
A semigroup $S$ is simple if and only if $S$ has an idempotent but
contains no isomorphic copy of the bicyclic semigroup $\C(p,q)$.
\end{theorem}

We recall that $\C(p,q)$ is a semigroup with a two-sided unit $1$,
generated by two elements $p,q$ and one relation $qp=1$.

Combining the Andersen Theorem~\ref{an}  with Theorem~\ref{t1} we
obtain another characterization of topological paragroups.

\begin{theorem}\label{t2}
A topological semigroup $S$ is a topological paragroup if and only
if $S$ is simple, contains an idempotent, contains no copy of the
bicyclic group, and each maximal subgroup $H_e$ of $S$ is a
topological group.
\end{theorem}

For compact topological semigroups the last three conditions
always are satisfied: such semigroups contain an idempotent by the
Iwassawa-Numakura Theorem~(see \cite{Iwassawa,
GelbaumKalishOlmstead, Numakura, Wallace1952} or \cite[Vol.~1,
Theorem~1.8]{CHK}), contain no copy of the bicyclic semigroup by
the Koch-Wallace Theorem~\cite{KochWallace} and the maximal
subgroups $H_e=eSe$ corresponding to minimal idempotents are
topological groups, being compact paratopological groups, see
\cite{Ellis}. In such a way we have proved the following theorem
due to Wallace \cite{Wallace}.

\begin{theorem}[Wallace]\label{t3}
Each simple compact topological semigroup $S$ is topologically
isomorphic to a topological paragroup.
\end{theorem}

In \cite[Theorem~I.5.3]{Rup} the Wallace Theorem was generalized
to compact semitopological semigroups.  By a {\em semitopological
semigroup} we understand a Hausdorff topological space $S$ endowed
with a separately continous semigroup operation.

\begin{theorem}[Ruppert]\label{t4}
Each simple compact semitopological semigroup $S$ is topologically
isomorphic to a topological paragroup.
\end{theorem}

Another direction of generalization of the Wallace Theorem consists in
replacing the compactness assumption by a weaker property. The
first step in this direction was made in \cite{GPR}.

\begin{theorem}[Gutik-Pagon-Repovs]\label{t5}
Each simple sequential countably compact topological semigroup is
a topological paragroup.
\end{theorem}

In this paper we generalize both the Wallace and
Gutik-Pagon-Repov\v{s} Theorems proving that simple topological
semigroups satisfying certain compactness-like properties are
topological paragroups. All topological spaces considered in this
paper are assumed to be Hausdorff.

We recall that a topological space $X$ is
\begin{itemize}
 \item {\em countably compact\/} if each closed discrete subspace of
 $X$ is finite;
 \item {\em pseudocompact\/} if $X$ is Tychonov and each continuous
 real-valued function on $X$ is bounded;
 \item {\em sequentially compact\/} if each sequence
 $\{x_n\}_{n\in\w}\subset X$
 has a convergent subsequence;
 \item {\em $p$-compact\/} for some free ultrafilter $p$ on $\w$ if
 each sequence $\{x_n\}_{n\in\w}\subset X$ has a $p$-limit
 $x_\infty=\lim_{n\to p}x_n$ in $X$.
\end{itemize}
Here the notation $x_\infty=\lim_{n\to p}x_n$ means that for each
neighborhood $O(x_\infty)\subset X$ of $x_\infty$ the set
$\{n\in\w:x_n\in O(x_\infty)\}$ belongs to the ultrafilter $p$. It
is clear that each sequentially compact and each compact
topological space is $p$-compact for every ultrafilter $p$.

By \cite{GS}, a topological space $X$ is $p$-compact for some free
ultrafilter $p$ on $\w$ if and only if each power $X^\kappa$ of
$X$ is countably compact if and only if the power $X^{2^\cc}$ is
countably compact. It is easy to see that each sequence $(x_n)_{n\in\w}$ in a countably compact topological space $X$ has $p$-limit $\lim_{n\to p}x_n$ for some free ultrafilter $p$ on $\w$. 

We shall say that for some free filter $p$ on $\w$ a double sequence $\{x_{m,n}\}_{m,n\in\w}\subset X$ has a double $p$-limit $\lim\limits_{n\to p}\lim\limits_{m\to p}x_{m,n}$ if $P=\{n\in\w:\exists \lim\limits_{m\to p}x_{m,n}\in X\}\in p$ and the sequence $(\lim\limits_{m\to p}x_{m,n})_{n\in P}$ has a $p$-limit in $X$.

We define a topological space $X$ to be {\em doubly countably compact} if each double sequence $(x_{m,n})_{m,n\in\w}$ in $X$ has a double $p$-limit $\lim\limits_{n\to p}\lim\limits_{m\to p}x_{m,n}\in X$ for some free ultrafilter $p$ on $\w$. 

\begin{proposition}\label{double} A topological space $X$ is doubly countably compact if $X$ is either sequentially compact or $p$-compact for some free ultrafilter $p$ on $\w$.
\end{proposition}

\begin{proof} The double countable compactness of $p$-compact spaces is obvious.
Now assume that $X$ is sequentially compact and take any double sequence $(x_{m,n})_{m,n\in\w}$  in $X$.  By the sequential compactness of $X$ there is an infinite subset $A_0\subset\w$ such that the subsequence $(x_{m,0})_{m\in A_0}$ converges to some point $x_0\in X$ in the sense that for each neighborhood $O(x_0)\subset X$ the set $\{n\in A_0:x_{m,0}\notin O(x_0)\}$ is finite. Now consider the sequence $(x_{m,1})_{m\in A_0}$ and by the sequential compactness of $X$ find an infinite subset $A_1\subset A_0$ such that the subsequence $(x_{m,1})_{m\in A_1}$ converges to some point $x_1$. Next we proceed by induction and for every $n\in\w$ construct an infinite subset $A_n\subset A_{n-1}$ such that the sequence $(x_{m,n})_{m\in A_n}$ converges to some point $x_n\in S$. Now take any infinite subset $A\subset \w$ such that for $A\subset^* A_n$ for every $n\in\w$. The latter means that the complement $A\setminus A_n$ is finite. It follows that for every $n\in\w$ the sequence $(x_{m,n})_{m\in A}$ converges to the point $x_n$. By the sequential compactness of $S$ for the sequence $(x_n)_{n\in A}$ there is an infinite subset $B\subset A$ such that the sequence $(x_n)_{n\in B}$ converges to some point $x\in X$.
Finally, take any free ultrafilter $p\ni B$ and observe that 
$x=\lim\limits_{n\to p}\lim\limits_{m\to p}x_{m,n}$.
\end{proof}

Theorem~\ref{t2} ensures that a simple topological semigroup $S$
is a topological paragroup provided
\begin{enumerate}
\item $S$ has an idempotent;
\item $S$ contains no copy of the bicyclic semigroup;
\item all maximal subgroups of $S$ are topological groups.
\end{enumerate}
\smallskip

Topological semigroups containing an idempotent can be characterized as follows.

\begin{theorem}\label{idempotent_char} A topological semigroup $S$ contains an idempotent if and only if for some $x\in S$  the double sequence $(x^{m-n})_{m\ge n}$ has a double $p$-limit $\lim\limits_{n\to p}\lim\limits_{m\to p}x^{m-n}\in S$ for some free ultrafilter $p$ on $\w$.
\end{theorem}  

\begin{proof} The ``only if'' part is trivial: just take any idempotent $x$ of $X$ and observe that $\lim\limits_{n\to p}\lim\limits_{m\to p}x_{m,n}=x$ for any free ultrafilter $p$ on $\w$.

To prove the ``if'' part, assume that for some $x\in S$ the double sequence $(x^{m-n})_{m\ge n}$ has a double $p$-limit $e=\lim\limits_{n\to p}\lim\limits_{m\to p}x^{m-n}$ for some free ultrafilter $p$ on $\w$. 

We claim that $e$ is an idempotent. Let $P\in p$ be the set of the numbers $n$ for which there is a $p$-limit $e_{-n}=\lim\limits_{m\to p}x^{m-n}$ in $S$. Then $e=\lim\limits_{P\ni n\to p}e_{-n}$.

Assuming that $e$ fails to be an idempotent, we can find a neighborhood $O(e)\subset S$ of $e$ such that $O(e)\cdot O(e)$ is disjoint with $O(e)$. Since $e=\lim\limits_{P\ni n\to p}e_{-n}$, the set $P_1=\{n\in P:e_{-n}\in O(e)\}$ belongs to the ultrafilter $p$.

Take any element $n\in P_1$ and observe that $\lim_{m\to p}x^{m-n}=e_{-n}\in O(e)$ implies $P_2=\{m\in P_1:m>n$ and $x^{m-n}\in O(e)\}\in p$. Pick any $m>n$ in $P_1$ and observe that $\lim_{i\to p}x^{i-m}=e_{-m}\in O(e)$ and thus the set $P_3=\{i\in P_2:i>m$ and $x^{i-m}\in O(e)\}$ belongs to $p$. Now take any number $i\in P_3$ and observe that 
$i\in P_3\subset P_2$ and $m\in P_2$ imply $x^{i-n},x^{i-m},x^{m-n}\in O(e)$. On the other hand, $x^{i-n}=x^{i-m}x^{m-n}\in O(e)\cdot O(e)\subset S\setminus O(e)$, which is a desired contradiction.
\end{proof} 

This characterization will be applied to obtain some convenient conditions on a topological semigroup $X$ guaranteeing the existence of an idempotent $e\in S$.

\begin{theorem}\label{idempotent_cor} A topological semigroup $S$ contains an idempotent if $S$ satisfies one of the following conditions:
\begin{enumerate}
\item $S$ is doubly countably compact;
\item $S$ is sequentially compact;
\item $S$ is $p$-compact for some free ultrafilter $p$ on $\w$;
\item $S^{2^{\mathfrak c}}$ is countably compact;
\item $S^{\kappa^\omega}$ is countably compact, where $\kappa$ is the minimal cardinality of a closed subsemigroup of $S$.
\end{enumerate}
\end{theorem}

\begin{proof} The first item follows immediately from Theorem~\ref{idempotent_char} and the definiton of a doubly sequentially countably compact space.

The next two assertions follow from the first one and Proposition~\ref{double}.
The fourth assertion follows from the third one and the characterization of spaces with countably compact power $S^{2^\mathfrak c}$ as $p$-compact spaces for some free ultrafilter $p$, see \cite{GS}. 

It remains to prove the last assertion. Let $\kappa$ be the smallest cardinality of a closed subsemigroup of $S$ and assume that the power $S^{\kappa^\w}$ is countably compact. Replacing $S$ by a suitable closed subsemigroup, we can assume that $|S|=\kappa$. Now it suffices to prove that the space $S$ is $p$-compact for some free ultrafilter $p$. For every $n\in\w$ consider the functional $\delta_n:S^\w\to S$ assigning to each function $f\in S^\w$ its value $\delta_n(f)=f(n)$ at $n$. This functional is an element of the power $S^{S^\w}$. The countable compactness of $S^{S^\w}$ guarantees that the sequence $(\delta_n)_{n\in\w}$ has an accumulation point $\delta_\infty\in S^{S^\w}$ and hence $\delta_\infty=\lim_{n\to p}\delta_n$ for some free ultrafilter $p$ on $\w$.
Then every function $f\in S^\w$ has the $p$-limit 
$$\lim_{n\to p}f(n)=\lim_{n\to p}\delta_n(f)=\delta_\infty(f),$$
which means that the space $S$ is $p$-compact.
\end{proof}

\begin{remark} Theorem~\ref{idempotent_cor} generalizes many known results related to idempotents in topological semigroups. In particular, it generalizes a result of A.~Tomita \cite{Tom96} on the existence of idempotents in $p$-compact cancellative semigroups as well as the classical Iwassawa-Numakura Theorem \cite[Vol.1, Theorem 1.6]{CHK} on the existence of an idempotent in compact topological semigroups.
\end{remark}  

The following theorem generalizing both the Wallace and
Gutik-Pagon-Repov\v{s} Theorems is the main result of this note.

\begin{theorem}\label{t12}
A simple topological semigroup $S$ is a topological paragroup if $S$ is doubly countably compact and has countably compact square $S\times S$.
\end{theorem}

This theorem follows immediately from Theorem~\ref{idempotent_char} and the
following characterizing

\begin{theorem}\label{t10}
A topological semigroup $S$ with countably compact square $S\times
S$ is a topological paragroup if and only if $S$ is simple and
contains an idempotent.
\end{theorem}

\begin{proof}
The ``only if'' part is trivial. To prove the ``if'' part, assume
that $S$ is simple and contains an idempotent.

First we check that $S$ contain no copy of $\C(p,q)$. Assume
conversely that $\C(p,q)\subset S$ and consider the sequence
$\{(q^n,p^n)_{n=1}^\infty\}$ in $\C(p,q)\times\C(p,q)\subset
S\times S$. The countable compactness of $S\times S$ guarantees
that this sequence has an accumulation point $(a,b)\in S\times S$.
Since $q^np^n=1$, the continuity of the semigroup operation on $S$
guarantees that $ab=1$. By Corollary~I.2~\cite{ES}, the bicyclic
semigroup $\C(p,q)$ endowed with the topology induced from $S$ is
a discrete topological space. So, we can find a neighborhood
$O(1)\subset S$ of $1\in\C(p,q)$ containing no other points of the
semigroup $\C(p,q)$. Since $ab=1$, the points $a,b$ have
neighborhoods $O(a),O(b)\subset S$ such that $O(a)\cdot
O(b)\subset O(1)$. Since $a$ is an accumulation point of the
sequence $q^n$, we can find $n\in\IN$ with $q^n\in O(a)$. By the
same reason, there is a number $m>n$ such that $p^m\in O(b)$. Then
$q^np^m=p^{m-n}\in O(a)\cdot O(b)\cap \C(p,q)=\{1\}$, which is a
contradiction. This contradiction shows that the simple semigroup
$S$ contains no copy of $\C(p,q)$ and hence is completely simple
by the Andersen's Theorem~\ref{an}.

For each idempotent $e\in E$ the maximal semigroup $H_e=eSe$ is
countably compact, being a continuous image of the countably
compact space $S$. Moreover, the square $H_e\times H_e$ is
countably compact, being a continuous image of the countably
compact space $S\times S$. Then $H_e$ is a topological group,
being a paratopological group with countably compact square, see
\cite{RR} or \cite[2.2]{AS}. Now Theorem~\ref{t2} assures that $S$
is a topological paragroup.
\end{proof}

For Tychonov topological semigroups the countable compactness of
the square $S\times S$ in the preceding theorem can be replaced by
its pseudocompactness.

\begin{theorem}
A topological semigroup $S$ with pseudocompact square $S\times S$
is a topological paragroup if and only if $S$ is simple and
contains an idempotent.
\end{theorem}

\begin{proof}
The ``only if'' part of the theorem is trivial. To prove the
``only if'' part, assume that the square $S\times S$ is
pseudocompact. By \cite[1.3]{BD}, the Stone-\v Cech
compactification $\beta S$ of $S$ is a compact topological
semigroup. By the Koch-Wallace Theorem \cite{KochWallace}, the
topological semigroup $\beta S$, being compact, contains no
isomorphic copy of the bicyclic semigroup and consequently, and so
does the subsemigroup $S$ of $\beta S$. By the Andersen Theorem
\cite[Theorem~2.54]{CP}, the simple semigroup $S$ is completely
simple. In order to apply Theorem~\ref{t1}, it remains to prove
that each maximal subgroup $H_e$ of $S$ is a topological group.
Since the idempotent $e$ of $S$ is primitive, the maximal group
$H_e$ coincides with $eSe$ and hence pseudocompact, being the
continuous image of the pseudocompact space $S$. Applying Theorem
2.6 of \cite{Rez}, we conclude that $H_e$, being a pseudocompact
paratopological group, is a topological group.
\end{proof}

For completely simple semigroups $S$ the pseudocompactness of the
square $S\times S$ in the preceding theorem can be replaced by the
pseudocompactness of $S$.

\begin{theorem}
A pseudocompact topological semigroup $S$ is a topological
paragroup if and only if $S$ is completely simple.
\end{theorem}

\begin{proof}
The ``only if'' part is trivial. To prove the ``if'' part, assume
that $S$ is a completely simple pseudocompact topological
semigroup. Take any primitive idempotent $e$ of $S$ and observe
that the maximal subgroup $H_e=eSe$ is pseudocompact, being the
continuous image of the pseudocompact space $S$. By Theorem 2.6 of
\cite{Rez}, the paratopological group $H_e$, being pseudocompact,
is a topological group. Applying Theorem~\ref{t1}, we conclude
that $S$ is a topological paragroup.
\end{proof}

Our final result describes the structure of simple sequential countably compact topological semigroups. We recall that a topological space $X$ is called {\em sequential} if for each non-closed subset $A\subset X$ there is a sequence $\{a_n\}_{n\in\w}\subset A$ that converges to some point $x\in X\setminus A$. 

\begin{theorem}\label{theorem16} For a simple topological semigroup $S$ the following conditions are equivalent:
\begin{enumerate}
\item $S$ is a regular sequential countably compact topological space;
\item $S$ is topologically isomorphic
to a topological paragroup $[X,G,Y]_{\sigma}$ for some
regular sequential countably compact topological spaces $X$ and $Y$
and a sequential countably compact topological group $G$.
\end{enumerate}
\end{theorem}

\begin{proof} $(1)\Ra(2)$. Assume that $S$ is a regular sequential countably compact topological space. It follows that $S$ is sequentially compact. By  Theorem~\ref{t12}, $S$ is topologically isomorphic to a topological paragroup $[X,G,Y]_{\sigma}$ for some
 topological spaces $X$ and $Y$
and some topological group $G$. The spaces $X,Y,G$, being homeomorphic to closed subspaces of $S$, are regular sequential and countably compact. 

$(2)\Ra(1)$ Assume that  $S$ is topologically isomorphic
to a topological paragroup $[X,G,Y]_{\sigma}$ for some
regular sequential countably compact topological spaces $X,Y$
and some sequential countably compact topological group $G$. It is clear that those spaces are sequentially compact. 
It follows from Boehme Theorem \cite[3.10.J(c)]{En} that the product $X\times G\times Y$ is sequential. Since the product of sequentially compact spaces is sequentially compact (and hence countably compact), the space $S$, being homeomorphic to $X\times G\times Y$, is regular sequential and countably compact.
\end{proof}


\section*{Open Problems}

\begin{problem} Let $X$ be a doubly countably compact space. 
Is the square $X\times X$ countably compact?
\end{problem}

\begin{problem}
Let $S$ be a (simple) semitopological semigroup with countably
compact power $S^{\cc}$. Has $S$ an idempotent?
\end{problem}

\begin{problem}
Assume that a simple Tychonov countably compact topological
semigroup $S$ contains an idempotent. Is $S$ completely simple?
Equivalently, is $S$ a topological paragroup?
\end{problem}

It is known that each (topological) semigroup embeds into a
simple (topological) semigroup, see \cite{Bruck},
\cite[Theorem~8.45]{CP} and \cite{Gutik}.

\begin{problem}
Is it true that each countably compact topological semigroup
embeds into a simple countably compact topological semigroup?
\end{problem}

\end{document}